\documentclass[11pt]{article}
\usepackage{latexsym}
\newtheorem{theorem}{Theorem}

\newtheorem{proposition}{Proposition}

\newcommand{\qed}{\ \hfill\mbox{$\Box$}\vspace{\baselineskip}}
\newenvironment{proof}{\noindent {\bf Proof:}}{{\qed}}

\begin{document}

\title{On braxtopes, a class of generalized simplices}

\author{Margaret M. Bayer\thanks{Supported in part by a grant from the
        University of Kansas General Research Fund}\\
        Department of Mathematics\\
        University of Kansas\\
        Lawrence KS 66045-7523 USA\\ \and
        Tibor Bisztriczky\thanks{Supported in part by a Natural Sciences and
        Engineering Research Council of Canada Discovery Grant}\\
        Department of Mathematics and Statistics\\
        University of Calgary\\
        Calgary, Alberta, T2N 1N4 Canada}

\date{July 13, 2006}

\maketitle

\begin{abstract}
In a $d$-simplex every facet is a $(d-1)$-simplex.
We consider as generalized simplices other combinatorial classes of
polytopes, all of whose facets are in the class.
Cubes and multiplexes are two such classes of generalized simplices.
In this paper we study a new class, braxtopes, which arise as the
faces of periodically-cyclic Gale polytopes.
We give a geometric construction for these polytopes and various
combinatorial properties.
\end{abstract}

\section{Introduction}
In the study of combinatorial properties of convex polytopes,
best understood are the simplicial ones.  
Among simplicial polytopes, cyclic polytopes have played an
important role.
They have the largest number of facets among all
polytopes with given dimension and number of vertices.
The combinatorial study of nonsimplicial polytopes is hampered by
the difficulty of generating classes with varied combinatorial
structure.
The simplicial (and their duals, the simple) polytopes are, in 
some sense, an extremal class of polytopes; we need other
extremal classes to better understand combinatorial parameters
associated with polytopes.

One approach is to find nonsimplicial analogues of cyclic polytopes.
In \cite{bisz-ord,bisz-pcg}, Bisztriczky introduced two such classes:
the {\em ordinary} polytopes (of odd dimensions) and the
{\em periodically-cyclic Gale} polytopes (of even dimensions).
The faces of the ordinary polytopes themselves form an interesting
class of polytopes, the {\em multiplexes} \cite{bisz-mult,bisz-multrig}.
Ordinary polytopes were studied further in 
\cite{bayer-ordflag,bayer-ordshell,bayer-bru-stew,dinh}.
The periodically-cyclic Gale polytopes have until now been less
studied.
In this paper we begin a study of periodically-cyclic Gale polytopes
by examining the polytopes that arise as their faces, and that are,
as are multiplexes, generalizations of simplices.

\section {Definitions}
Let $Y$ be a set of points in ${\bf R}^d$, $d\ge 1$.
Then $[Y]$ and $\langle Y\rangle$ denote, respectively, the convex hull
and the affine hull of $Y$.
If $Y=\{y_1,y_2,\ldots, y_s\}$ is a finite set, we set 
$[y_1,y_2,\ldots,y_s]=[Y]$ and 
$\langle y_1,y_2,\ldots,y_s\rangle=\langle Y\rangle$.

Let $V=\{x_0,x_1,\ldots, x_n\}$ be a totally ordered set of $n+1$ points in 
${\bf R}^d$ with $x_i<x_j$ if and only if $i<j$.
We say that $x_j$ {\em separates} $x_i$ and $x_k$ 
if $x_i<x_j<x_k$.
For $Y\subset V$,
$Y$ is a {\em Gale set} (in $V$) if any two points of $V\setminus Y$ are
separated by an even number of points of $Y$.

Let $P\subset{\bf R}^d$ be a (convex) $d$-polytope.
For $-1\le i\le d$, let ${\cal F}_i(P)$ denote the set of $i$-dimensional
faces of $P$ and $f_i(P)=|{\cal F}_i(P)|$.
For convenience, we set ${\cal V}(P)={\cal F}_0(P)$, 
${\cal E}(P)={\cal F}_1(P)$, and ${\cal F}(P)={\cal F}_{d-1}(P)$.
We recall that $(f_{-1}(P),f_0(P),f_1(P),\ldots, f_{d-1}(P), f_d(P))$ is the
{\em $f$-vector} of $P$.
In the case that $P$ is simplicial, the {\em $h$-vector} of $P$ is
$(h_0(P),h_1(P),\ldots, h_d(P))$ with 
$h_i(P)=\sum_{j=0}^i (-1)^{i-j}{d-j\choose d-i}f_{j-1}(P)$.
A chain of faces
$\emptyset\subset G_1\subset G_2\subset \cdots \subset G_r\subset P$
is an {\em $S$-flag} if $S=\{\dim G_1,\dim G_2,\ldots, \dim G_r\}$.
Writing $f_S(P)$ as the number of $S$-flags of $P$, the {\em flag vector}
$(f_S(P))_{S\subseteq\{0,1,\ldots, d-1\}}$ of $P$ is a vector with $2^d$
entries.
Finally we refer to \cite{ziegler} for the definitions of a {\em triangulation}
and a {\em shelling} of $P$.

Let ${\cal V}(P)=\{x_0,x_1,\ldots, x_n\}$, $n\ge d$.
We set $x_i<x_j$ if and only if $i<j$, and call $x_0<x_1<\cdots <x_n$ a
{\em vertex array} of $P$.
Let $G\in{\cal F}_i(P)$, $1\le i\le d-1$, such that 
$G\cap{\cal V}(P)=\{y_0,y_1,\ldots, y_m\}$ (each $y_j$ is some $x_i$).
Then $y_0<y_1<\cdots <y_m$  is the vertex array of $G$ if it is induced by
$x_0<x_1<\cdots <x_n$. 
Finally, $P$ with $x_0<x_1<\cdots <x_n$  is {\em Gale} (with respect to the
vertex array) if the vertex set of each facet of $P$ is a Gale set.

We recall from \cite{gale} and \cite{grunbaum} that a $d$-polytope $P$ is 
{\em cyclic} if it is simplicial and Gale with respect to some vertex array.
From \cite{bisz-pcg}, $P$ is {\em periodically-cyclic} if there is a vertex
array, say, $x_0<x_1<\cdots <x_n$ and an integer $k$ with $n+1\ge k\ge d+2$,
such that
\begin{itemize}
\item $[x_{i+1},x_{i+2},\ldots, x_{i+k}]$ is a cyclic $d$-polytope with
      $x_{i+1}<x_{i+2}<\cdots<x_{i+k}$, for $-1\le i\le n-k$,  and
\item $[x_{i+1},x_{i+2},\ldots, x_{i+k+1}]$ is not cyclic 
      for $-1\le i\le n-k-1$.
\end{itemize}
We note that if $k=n+1$, then $P$ is cyclic.

From \cite{bisz-mult}, $P$ is a {\em multiplex} if there is a vertex
array, say, $x_0<x_1<\cdots <x_n$ such that the facets of $P$ are 
$[x_{i-d+1},\ldots, x_{i-1},x_{i+1},\ldots, x_{i+d-1}]$ for 
$0\le i\le n$ under the convention: $x_t=x_0$ for $t\le 0$ and
$x_t=x_n$ for $t\ge n$.
Next, $P$ is {\em multiplicial} if each facet of $P$ is a multiplex with
respect to the ordering induced by a fixed vertex array of $P$.
Finally, $P$ is {\em ordinary} if it is Gale and multiplicial with respect
to some vertex array.
We note from \cite{bisz-ord} that if $d\ge 4$ is even, then every ordinary 
$d$-polytope is cyclic.

We observe that the noncyclic polytopes mentioned above are nonsimplicial.
In this spirit we introduce another class of nonsimplicial polytopes, the
braxtopes.

\vspace{1\baselineskip}

\noindent {\bf Definition.}
For $d\le 2$, a {\em $d$-braxtope} is a $d$-simplex.

For $n\ge d\ge 3$, $P$ is a  {\em $d$-braxtope} if there is a
vertex array, say, $x_0 < x_1 < \ldots < x_n$ 
such that the facets of $P$ are 
$$T_i=[x_i,x_{i+1},\ldots, x_{i+d-1}] \mbox{ for $0\le i\le n-d+1$}$$
and
$$E_j = [x_0,x_{j-(d-2)},\ldots,x_{j-1},x_{j+1}, \ldots, x_{j+(d-2)}] 
\mbox{ for $2\le j\le n$}$$
under the convention: $x_t=x_0$ for $t\le 0$ and $x_t=x_n$ for $t\ge n$.

\vspace{1\baselineskip}

We note that 
$E_n=[x_0,x_{n-(d-2)},\ldots, x_{n-1},x_n]$ and $|{\cal F}(P)|=2n-d+1$.
Finally, $P$ is {\em braxial} if each proper face of $P$ is a braxtope with 
respect to the ordering induced by a fixed vertex array of $P$.

\section{Realizability and properties of braxtopes}
Henceforth $Q^{d,n}$ denotes a $d$-braxtope with the vertex array
$x_0<x_1<\cdots <x_n$, $n\ge d\ge 3$.
For $n=d$, $Q^{d,d}$ is a $d$-simplex.
If  $d+1\le n\le 2d-2$, then $Q^{d,n}$ is a face
of a periodically-cyclic Gale $2m$-polytope when $2m\ge d+1$ \cite{bisz-pcg}.  
This is not obvious, and this observation by the first author led to the
formulation of braxtopes as a new class of polytopes.

\begin{theorem}
$Q^{d,n}$ is realizable in ${\bf R}^d$ for all $n\ge d\ge 3$.
\end{theorem}

\begin{proof}
In view of the preceding, we may assume that $n\ge 2d-1$ and that 
$Q^{d,n-1}\subset{\bf R}^d$ exists with $x_0<x_1<\cdots <x_{n-1}$.
Let $Q'=Q^{d,n-1}$ and 
${\cal F}(Q')=\{T'_0,\ldots, T'_{n-d},E'_2,\ldots, E'_{n-1}\}$.
It is easy to check by a simple beneath-beyond argument (see \cite{grunbaum})
that $[Q',x_n]$ is a $Q^{d,n}$ with $x_0<x_1<\cdots <x_{n-1}<x_n$ if
$x_n\in{\bf R}^d$ is a point with the properties:
\begin{itemize}
\item $\displaystyle x_n\in L=\langle x_0,x_{n-d},x_{n-d+1},x_{n-1}\rangle=
      \bigcap_{j=n-d+2}^{n-2} E'_j$,
\item $x_n$ is beyond $E'_{n-1}$, and
\item $x_n$ is beneath every other facet of $Q'$.
\end{itemize}
Specifically, the three $(d-2)$-faces $[x_{n-d+1},\ldots, x_{n-1}]$,
$[x_0,x_{n-d+1},\ldots, x_{n-2}]$, and 
$[x_0,x_{n-d+2},\ldots, x_{n-1}]$ of $E'_{n-1}$ yield $T_{n-d+1}$, $E_{n-1}$,
and $E_n$, respectively.
We observe that the existence of such a point
$x_n$ is due to the fact that $L$ is a 3-flat,
$L\cap\langle E'_{n-1}\rangle=\langle x_0,x_{n-d+1},x_{n-1}\rangle$ is a supporting
plane of $L\cap Q'$, and 
$|{\cal F}(Q')\setminus\{E'_{n-d+2},\ldots,E'_{n-1}\}|$ is finite.
\end{proof}

\begin{proposition}\label{prop1}
Let $Q=Q^{d,n}$ be a $d$-braxtope with $x_0<x_1<\cdots <x_n$, $n\ge d\ge 3$.
Then 
\begin{enumerate}
\item $[x_0,x_u]\in{\cal E}(Q)$ for $1\le u\le n$,
      \label{edges_on_0}
\item $[x_1,x_u]\in{\cal E}(Q)$ if and only if $u=0$ or $2\le u\le d$,
      \label{edges_on_1}
\item $[x_u,x_n]\in{\cal E}(Q)$ if and only if $u=0$ or $n-(d-1)\le u\le n-1$,
      \label{edges_on_n}
\item for $2\le t\le n-1$, $[x_t,x_u]\in{\cal E}(Q)$ if and only if $u=0$, or
      $t-d+1\le u\le t+d-1$, $u\ne t$,\label{middle_edges}
\item $[x_0,x_{t+1},x_{t+k}]\in{\cal F}_2(Q)$ for $0\le t\le n-k$ 
      and $2\le k\le d-2$,\label{2-faces}
\item $[x_0,x_t,x_{t+1},x_{t+d-1},x_{t+d}]\in{\cal F}_3(Q)$ for $1\le t\le n-d$,
      and\label{3-faces}
\item $\{x_t,x_{t+1},\ldots x_{t+d}\}$ is an affinely independent set for
      $0\le t\le n-d$.\label{affind}
\end{enumerate}
\end{proposition}

\begin{proof}
(\ref{edges_on_0}) We check that each edge
    $[x_0,x_u]$ is the intersection of specific facets.
For example, $\displaystyle[x_0,x_1]=T_0\cap \bigcap_{j=2}^{d-1} E_j$ and 
$[x_0,x_u]=E_{u-d+2}\cap E_{u+d-2}$ for $d\le u\le n-d+1$.

(\ref{edges_on_1}) We note $T_1=[x_1,x_2,\ldots, x_d]\in{\cal F}(Q)$ and for 
$u\ge d+1$, at most $d-3$ facets of $Q$ contain $[x_1,x_u]$.

(\ref{edges_on_n}) and (\ref{middle_edges}) We argue as in (\ref{edges_on_1}).

(\ref{2-faces}) If $t+1\ge n-d+2$ or $t+k\le d-1$, then $[x_0,x_{t+1},x_{t+k}]$
is in a simplex facet ($E_n$ or $T_0$), and so it is a 2-face.
Assume $d-k\le t\le n-d$.  Then 
\begin{eqnarray*}
\lefteqn{E_{t-d+k+2}\cap E_{t+d-1}\cap \bigcap_{j=t+2}^{t+k-1} E_j}\\
      &=& [x_0,x_{t+1},x_{t+2},\ldots, x_{t+k}]\cap 
          \bigcap_{j=t+2}^{t+k-1} E_j \\
      &=& [x_0,x_{t+1},x_{t+k}]
\end{eqnarray*}
So $ [x_0,x_{t+1},x_{t+k}]\in {\cal F}_2(Q)$.

(\ref{3-faces}) $\displaystyle[x_0,x_t,x_{t+1},x_{t+d-1},x_{t+d}]
=\bigcap_{j=t+2}^{t+d-2}E_j$,
$[x_t,x_{t+d}]\not\in{\cal E}(P)$, and $[x_0,x_t,x_{t+1}]\in{\cal F}_2(P)$.

(\ref{affind}) This follows immediately from 
    $\{T_0,T_1,\ldots, T_{n-d}\}\subset{\cal F}(Q)$.
\end{proof}

\begin{theorem}
Let $Q=Q^{d,n}$ with 
$x_0 < x_1 < \ldots < x_n$, $n\ge d+1\ge 4$. 
Then $Q'=[x_0,x_1,\ldots, x_{n-1}]$ is a $d$-braxtope with 
$x_0 < x_1 < \ldots < x_{n-1}$.
\end{theorem}

\begin{proof}
With the notation above, we observe that 
$$\{T_0,\ldots, T_{n-d-1},E_2,\ldots, E_{n-d+1}\}\subset {\cal F}(Q').$$
Let $n-d+2\le j\le n-2$.
Then $x_n\in E_j$ yields that 
$$E'_j=[x_0,x_{j-d+2},\ldots, x_{j-1},x_{j+1},\ldots, x_{n-1}]\in{\cal F}(Q').$$
Thus, we need only to verify that 
$E'_{n-1}=[x_0,x_{n-d+1},\ldots, x_{n-1}]\in{\cal F}(Q')$.
By Proposition~\ref{prop1}(\ref{edges_on_n}), $x_n$ is a simple vertex of $Q$ and 
$[x_u,x_n]\in{\cal E}(Q)$ for exactly 
$x_u\in X=\{x_0,x_{n-d+1},\ldots, x_{n-1}\}$.
By Proposition~\ref{prop1}(\ref{affind}), $\langle X\rangle$ is a hyperplane
of ${\bf R}^d$ and $[X]$ is a $(d-1)$-polytope.

Finally, $x_n\not\in Q'$ implies that $x_n$ is beyond some $F'\in{\cal F}(Q')$.
Let $x_u$ be a vertex of $F'$.
Since there is an $F\in{\cal F}(Q)$ such that $x_u\in F$ and $x_n\not\in F$,
it follows that $[x_u,x_n]\in{\cal E}(P)$ and $x_u\in X$.
Thus $|X|=d$ yields that $F'=[X]=E'_{n-1}$, and ${\cal F}(Q')$ contains the
set of facets of a $Q^{d,n-1}$.
It is well known that this implies $Q'=Q^{d,n-1}$.
\end{proof}

\begin{theorem} \label{braxial} 
Let $Q=Q^{d,n}$ with $x_0<x_1<\cdots <x_n$, $n\ge d\ge 3$.  Then
\begin{enumerate}
\item $Q$ is a braxial $d$-polytope.\label{Qbraxial}
\item The vertex figure $Q/x_0$ of $Q$ at $x_0$ is a $(d-1)$-multiplex
      with the induced ordering.\label{verfig_0}
\item Let $n\le 2d-3$. Then $Q$ is a $(2d-2-n)$-fold pyramid
      over an $(n-d+2)$-braxtope with the induced ordering.\label{pyramid}
\item For $-1\le j\le d$, $f_j(Q)={d+1\choose j+1}+(n-d)\left[{d-1\choose j}
      + {d-2\choose j-1}\right]$.\label{f-vector}
\item $Q$ is elementary; that is, $f_{\{0,2\}}(Q)-3f_2(Q)+f_1(Q)-df_0(Q)+
      {d+1\choose 2}=0$.\label{elementary} 
\end{enumerate}
\end{theorem}

\begin{proof}
(\ref{Qbraxial}) We verify that each $F\in{\cal F}(Q)$ is a $(d-1)$-braxtope
with the induced ordering.
Assume $F$ is not a $(d-1)$-simplex, and hence, 
$F\in\{E_3,\ldots, E_{n-2}\}$ with the standard notation.
If $f_0(F)=m+1$ and $F$ is a $(d-1)$-braxtope, then we
denote its $(d-2)$-faces by $T'_0,\ldots, T'_{m-d+2},E'_2,\ldots, E'_m$.

If $3\le u\le d-1$, then 
$E_u=[x_0,x_1,\ldots,x_{u-1},x_{u+1},\ldots,x_{u+d-2}]$,
$m=u+d-3$, $T'_i=T_i\cap E_u$ for $0\le i\le u-1=m-d+2$, 
$\{E'_2,E'_3,\ldots, E'_{u+d-4}\}=\{E_j\cap E_u|\mbox{ $2\le j\le u+d-3$,
$j\ne u$ }\}$, and $E'_m=E_u\cap E_{u+d-1}$.
A similar argument yields the claim for $E_d$, $E_{d+1}$, \ldots, $E_{n-2}$.

(\ref{verfig_0}) Since $[x_0,x_u]\in{\cal E}(P)$ for $1\le u\le n$, it follows
from the description of ${\cal F}(Q)$ that the $(d-2)$-faces of
$\overline{Q}=Q/x_0$ are (writing $\overline{x}_u$ for the vertex of 
$\overline{Q}$ corresponding to $[x_0,x_u]$),
$[\overline{x}_{i-d+2},\ldots,\overline{x}_{i-1},\overline{x}_{i+1},\ldots,
\overline{x}_{i+d-2}]$ for $1\le i\le n$, with the convention 
$\overline{x}_r=\overline{x}_1$ for $r\le 1$ and 
$\overline{x}_r=\overline{x}_n$ for $r\ge n$.
These are the $(d-2)$-faces of a $(d-1)$-multiplex with $n$ vertices.

(\ref{pyramid}) We observe that for $n\le 2d-3$, $Q=[E_{n-d+2},x_{n-d+2}]$; that
is, $Q$ is a pyramid over the $(d-1)$-braxtope $E_{n-d+2}$ with apex 
$x_{n-d+2}$.
We note that $E_{n-d+2}=Q^{d-1,n-1}$ with $n-1\le 2d-4=2(d-1)-2$.
Thus, either $n=2d-3$ and we are done, or $n-1\le 2(d-1)-3$, and we
repeat the argument.
In summary, $Q$ is a $(2d-2-n)$-fold pyramid over the $(n-d+2)$-braxtope
$[x_0,x_1,\ldots, x_{n-d+1},x_d,\ldots, x_n]$ with apices
$x_{n-d+2},\ldots, x_{d-2}, x_{d-1}$.

(\ref{f-vector}) 
We count the faces of $Q^{d,n}$ in two groups: those faces containing
the vertex $x_0$, and those not containing $x_0$.
The former are intersections of the facets $E_i$ (and $T_0$), and the latter 
are all contained in the facets $T_i$ ($1\le i\le n-d+1$).

Recall that the vertex figure of $x_0$ in $Q$ is the 
$(d-1)$-multiplex with $n$ vertices.
The $f$-vector of the multiplex is given in \cite{bisz-mult}.
Thus the number of $j$-faces of $Q$ containing $x_0$ is 
${d\choose j}+(n-d){d-2\choose j-1}$.

For $d+1\le \ell\le n$, ${d-1\choose j}$ is the number of $j$-faces of
$\displaystyle\bigcup_{i=1}^{n-d+1} T_i$ containing $x_\ell$ as the greatest
vertex.
The number of $j$-faces in $T_1=[x_1,x_2,\ldots, x_d]$ is 
${d\choose j+1}$.
Thus the total number of $j$-faces in 
$\displaystyle\bigcup_{i=1}^{n-d+1} T_i$
is ${d\choose j+1}+(n-d){d-1\choose j}$.

(\ref{elementary}) 
In the inductive construction of the $d$-braxtope (proof of Theorem~A),
new vertices are not placed on flats spanned by 2-faces.
So all 2-dimensional faces of every $d$-braxtope are triangles, and thus
$f_{\{0,2\}}(Q)-3f_2(Q)=0$.
Now $$f_1(Q)-df_0(Q)+{d+1\choose 2}={d+1\choose 2}+(n-d)d-d(n+1)+{d+1\choose 2}=0.$$
\end{proof}

\noindent {\bf Remarks.}
If $n>d$, then the $f$-vector of the $d$-braxtope equals the $f$-vector of the
$(d-3)$-fold pyramid over the bipyramid over an $(n-d+2)$-gon.
We conjecture that this result extends to flag vectors.

Kalai \cite{Kalai-NATO} introduced {\em elementary} polytopes as $d$-polytopes 
satisfying $f_{\{0,2\}}-3f_2+f_1-df_0+ {d+1\choose 2}=0$.  
This quantity is nonnegative for all $d$-polytopes by a rigidity argument
\cite{kalai-rigid}.
It may be interpreted as the difference $h_2-h_1$ of (middle perversity) betti 
numbers of the associated toric variety.

\section{Triangulation and the $h$-vector of the braxtope}

Earlier we defined the $h$-vector of a simplicial polytope by a linear
transformation of the $f$-vector.
The definitions of $f$-vector and $h$-vector extend naturally to 
simplicial complexes.
The $h$-vector of a simplicial polytope is the sequence of cohomology
ranks of the toric variety associated to the polytope.
For nonsimplicial polytopes the middle perversity betti numbers of the
toric variety form the $h$-vector, but it depends on the flag vector,
not just on the $f$-vector.
A triangulation $\Delta$ of a polytope $P$ is {\em shallow} if every
$k$-face of $\Delta$ is contained in a face of $P$ of dimension
at most $2k$.
The $h$-vector of a shallow triangulation (if one exists) may be used to 
compute the $h$-vector of the nonsimplicial polytope \cite{bayer-weakly}.

\begin{theorem}
Let $Q=Q^{d,n}$ with 
$x_0 < x_1 < \ldots < x_n$, $n\ge d\ge 3$. 
\begin{enumerate}
\item For $1\le i\le n-d+1$, let
      $$J_i=[x_0,x_i,x_{i+1},\ldots, x_{i+d-1}].$$
      Then $\{J_1,J_2,\ldots, J_{n-d+1}\}$ are the facets of a triangulation 
      $\Delta$ of $Q$.\label{triangulation}
\item $\Delta$ is a shallow triangulation of $Q$.\label{shallow}
\item $h(Q)=(1,n-d+1,n-d+1,\ldots, n-d+1,1)$.\label{h-vector}
\end{enumerate}
\end{theorem}
\begin{proof}
(\ref{triangulation}) Since the facets not containing $x_0$ are the simplices
$T_i$ ($1\le i\le n-d+1$),
this is the triangulation resulting from {\em pulling}
the vertex $x_0$ (see \cite{lee-triang}).  

(\ref{shallow}) 
First observe that for any $j$, $2\le j\le n-d+2$,
$\displaystyle \bigcap_{i=j}^{j+d-3}E_i=[x_0,x_{j-1},x_{j+d-2}]$.
In particular this intersection is two-dimensional, and so for any
set $I$ contained in a consecutive $(d-2)$-element subset of 
$\{2,3,\ldots, n\}$, $\displaystyle \dim \bigcap_{i\in I}E_i= d-|I|$. 

Suppose $\sigma$ is a face of $\Delta$.
Note that all vertices and edges of $\Delta$ are vertices and edges of
$Q$, since $[x_0,x_i]\in{\cal E}(Q)$ for all $i$, so assume
$\dim\sigma\ge 2$.
If $x_0\not\in\sigma$, then $\sigma\subset T_i$ for some facet $T_i$  ($i\ge 1$)
of $Q$, and $\sigma$ is thus a face of $Q$.
Now suppose $x_0\in \sigma\subset[T_{i+1},x_0]$.
Let $\tau=[\sigma,x_{i+1},x_{i+d}]$, and $I=\{j:\mbox{$i+1\le j\le i+d$
, $x_j\not\in \tau$}\}$.
Then $\displaystyle\sigma\subseteq\tau\subseteq\bigcap_{i\in I}E_i$, and
$$\dim\bigcap_{i\in I}E_i=d-|I|=f_0(\tau)-1
\le  f_0(\sigma)+1=\dim\sigma+2\le 2\dim\sigma.$$
Thus $\sigma$ is contained in the face 
$\displaystyle\bigcap_{i\in I}E_i$ of $Q$ of dimension at most $2\dim\sigma$.

(\ref{h-vector}) The ordering $J_1,J_2,\ldots, J_{n-d+1}$ of the facets of
$\Delta$ forms a shelling of $\Delta$; for $2\le j\le n-d+1$, the unique
minimal face of $\displaystyle J_j\setminus\bigcup_{i<j} J_i$ is
$\{x_{j+d-1}\}$.
So the $h$-vector of $\Delta$ is $(1,n-d,0,0,\ldots,0)$.
By \cite[Theorem~4]{bayer-weakly} this implies that the $h$-vector of $Q$
is $h(Q)=(1,n-d+1,n-d+1,\ldots, n-d+1,1)$.
\end{proof}

\noindent {\bf Remark.} The formula for $h(Q^{d,n})$ would also follow
from the conjecture that the flag vector of $Q^{d,n}$ equals the flag
vector of the $(d-3)$-fold pyramid over the bipyramid over the 
$(n-d+2)$-gon.

The colex order of the facets of $Q=Q^{d,n}$ gives a shelling of $Q$.
Like the colex shelling of the ordinary polytopes and 
multiplexes \cite{bayer-ordshell}, 
this shelling of $Q$ has special properties that are important for
counting faces: for each $j$, $F_j\setminus\cup_{i<j}F_i$ has a
unique minimal face $G_j$, which is a simplex, and the quotient
polytope $F_j/G_j$ is a simplex.

\section{Extension}
Theorem \ref{braxial}(\ref{verfig_0}) says that the vertex figure of $x_0$ in 
a braxtope $Q$ is a multiplex.
The antistar of $x_0$ (the polytopal complex of faces of $Q$ not
containing $x_0$) is a triangulation of the multiplex into the
simplices $T_1$, $T_2$, \ldots, $T_{n-d}$.
This multiplex-braxtope relationship may be extended by adding more
vertices like $x_0$.

\vspace{1\baselineskip}

\noindent {\bf Definition.}
For $d\le r+1$, an {\em $(r,d)$-braxtope} is a $d$-simplex.

For $n\ge d\ge r+2\ge 2$, $P$ is an  {\em $(r,d)$-braxtope} if there is a
vertex array, say, $x_0 < x_1 < \ldots < x_n$
such that the facets of $P$ are
$$T_{i,j}=[\{x_0,x_1,\ldots, x_{r-1}\}\setminus\{x_i\}\cup
\{x_j,x_{j+1},\ldots, x_{j+d-r}\}]$$
for $0\le i\le r-1\le r \le j\le n-d+r$,
$$T_{0,0}=[x_0,x_1,\ldots, x_{d-1}],$$
and
$$E_j = [\{x_0,x_1,\ldots, x_{r-1},x_{j-(d-r-1)},\ldots,x_{j-1},x_{j+1}, 
\ldots, x_{j+(d-r-1)}\}]$$ 
for $r+1\le j\le n$,
under the convention: $x_t=x_0$ for $t\le 0$ and $x_t=x_n$ for $t\ge n$.

\vspace{1\baselineskip}

We note that a $(1,d)$-braxtope is a $d$-braxtope.
If $r=0$, we understand that there are no facets $T_{i,j}$ (except $T_{0,0}$)
and that the set $\{x_0,x_1,\ldots, x_{r-1}\}$ is empty, so that a
$(0,d)$-braxtope is a $d$-multiplex.
The theorems in this paper have natural analogues for $(r,d)$-braxtopes.
The $(1,d)$-braxtopes are of special interest because they arise as facets
of periodically-cyclic Gale polytopes.
It would be interesting to investigate polytopes, all of whose facets are
$(r,d)$-braxtopes.

\end{document}